\documentclass[12pt]{amsart}
\usepackage{amsfonts}

\newtheorem{theorem}{Theorem}
\newtheorem{lemma}{Lemma}

\theoremstyle{definition}

\theoremstyle{remark}
\newtheorem{remark}[theorem]{Remark}

\begin{document}
\title[OPERATORS IN THE GRAND LEBESQUE SPACES]
{TENSOR, SOBOLEV, MULTIPLICATIVE AND CONVOLUTION\\ OPERATORS IN
          THE  BIDE - SIDE GRAND LEBESQUE SPACES}

\author{E. Liflyand}

\address{Department of Mathematics,
Bar-Ilan University, 52900 Ramat-Gan, Israel}
\email{liflyand@math.biu.ac.il}

\author{E. Ostrovsky}
\email{galo@list.ru}

\author{L. Sirota}
\email{sirota@zahav.net.il}

\begin{abstract}
   In this paper we study the multiplicative, tensor, Sobolev's and
 convolution inequalities in certain Banach spaces, the so-called
 Bide - Side Grand Lebesque Spaces, and give examples to show their
 sharpness.
\end{abstract}

\maketitle


\quad


\section{Introduction}

Let $ (X, \Sigma,\mu) $ be a measurable space with non - trivial
measure $\mu.$ We suppose the measure to be diffuse, that is, for
all $A \in \Sigma$  such that $\mu(A) \in (0, \infty)$ and there
exists $B \subset A$ such that $\mu(B) = \mu(A)/2.$ We also
suppose that the measure is $\sigma$-finite: there is a sequence
$E(n) \in \Sigma$ such that $\mu(E(n)) < \infty$ and
$\cup_{n=1}^{\infty}E(n) = X.$

For $a$ and $b$ constants, $1 \le a < b \le \infty,$ let $\psi =
\psi(p),$ $p \in (a,b)$ be a continuous log-convex positive
function such that $\psi(a + 0)$ and $\psi(b-0)$ exist, with
$\max(\psi(a+0), \psi(b-0)) = \infty$ and $\min(\psi(a+0),
\psi(b-0))> 0.$

The Bide-Side Grand Lebesque Space $ BSGLS(\psi; a,b)= G_X(\psi;
a,b) = G(\psi; a,b) = G(\psi) $ consists, by the well-known
definition, of all measurable functions $ h: X \to{\mathbb R} $
with finite norm

\begin{eqnarray*}
||h||_{G(\psi)} \stackrel{def}{=}\sup_{p \in (a,b)}|h|_p/\psi(p),
 \quad |h|_p = \left[\int_X |h(x)|^p \ \mu(dx) \right]^{1/p}.
\end{eqnarray*}
These spaces are intensively studied, in particular, their
associate spaces, fundamental functions $\phi(\delta; G(\psi;
a,b)),$ Fourier and singular operators, condition for convergence
and compactness, reflexivity and separability, martingales in
these spaces, etc.; see, e.g., \cite{BS1,CM3,
Fio4,FK5,IS7,IKO8,JM9,KM10,KPS12,OS14}. They are also Banach and
moreover rearrangement invariant (r.i.) spaces.

The BSGLS norm, in particular, Orlicz norm estimates for
measurable functions, e.g., for random variables are used in PDE
\cite{Fio4,FK5,IS7,IKO8}, probability in Banach spaces
\cite{LedTal18}, in the modern non-parametrical statistics, for
example, in the so-called regression problem \cite{Os13, OS14}.
The latter reads as follows. Given the observation $ \{\xi(i)\},$
$i = 1,2,3,\ldots,n,$ with $n \to \infty,$ of the view

\begin{eqnarray*}
\xi(i) = g(z(i)) + \epsilon(i), \ i = 1,2,\ldots,
\end{eqnarray*}
where $ g(\cdot) $ is an unknown estimated function, $
\{\epsilon(i) \} $ are the errors of measurements that may be
independent random variables or martingale differences,  $
\{z(i)\} $ is a dense set in the metric space $ (Z,\rho) $ with
Borel measure $ \nu. $ \par

 Let $ \{\phi_k(z) \} $ be a complete orthonormal sequence of functions,
for example, the classical trigonometric sequence, Legengre or
Hermite polynomials, etc. Put

\begin{eqnarray*}
c_k(n) = n^{-1} \sum_{i=1}^n \phi_k(z(i))&,& \ \tau(N) = \tau(N,n)
= \sum_{k = N+1}^{2N} (c_k(n))^2,\\
 M = \arg\min_{N \in [1,n/3]} \tau(N)&,& \ f_n(z) = \sum_{k=1}^M
c_k(n) \phi_k(z).                                \end{eqnarray*}
To study the confidence region for estimating function $ f $ in
the $ L(p)$ norm, written $ |f_n - f|_p,$ exponential bounds for
the tail of the distribution of polynomial martingales are used
being derived via the BSGLS spaces.
\par

 Let $ a \ge 1,$ $ b\in (a,\infty],$ and let $ \psi =
\psi(p) $ be a positive continuous function on the {\it open}
interval $ (a,b) $ such that  there exists a measurable function
$f: X \to{\mathbb R} $ for which

\begin{eqnarray}\label{no0}
f(\cdot) \in \cap_{p \in (a,b) } L_p, \ \psi(p) = |f|_p, \ p \in
(a,b).\end{eqnarray}
  We say that the equality (\ref{no0}) and the function $  f(\cdot) $ from
(\ref{no0}) is the {\it representation} of the function $ \psi. $

We denote the set of all these functions by $ \Psi: \ \Psi =
\Psi(a,b) = \{ \psi(\cdot) \}. $  For complete description of
these functions see, for example, (\cite[p.p. 21 - 27]{Os13},
\cite{OS14}).

\begin{remark} Observe that if $ \psi_1 \in \Psi(a_1,b_1)$ and
$\psi_2\in \Psi(a_2,b_2),$ with $(a_1,b_1)\cap (a_2,b_2) =
(a_3,b_3) \ne \emptyset,$ then $ \psi_1\psi_2\in \Psi(a_3,b_3).$
Indeed, if $\psi_1(p) = |f_1|_p$ and $\psi_2(p) = |f_2|_p,$ and
the functions $ f_1$ and $f_2 $ are independent in the
probabilistic sense, that is, for all Borel sets $ A,B $ on the
real axis $\mathbb R$

\begin{eqnarray*}
\mu \{x: \ f_1(x) \in A, \ f_2(x) \in B \} = \mu  \{x: f_1(x) \in
A \} \ \mu \{x: \ f_2(x) \in B \},
\end{eqnarray*}
then $\psi_1(p) \psi_2(p) = |f_1  f_2|_p.$ \end{remark}

We note that the $ G(\psi) $ spaces are also interpolation spaces
(the so-called  $\Sigma$-spaces), see
\cite{BS1,CM3,Fio4,FK5,IS7,IKO8,JM9,KPS12,OS14}, etc. However, we
hope that our direct representation of these spaces is of certain
convenience in both theory and applications.

The  $G(\psi)$ spaces with $ \mu(X) = 1$ appeared in \cite{KO11};
it was proved that in this case each $ G(\psi) $ space coincides
with certain exponential Orlicz space, up to norm equivalence.
\par

{\bf  The main goal of this paper is to prove new (and extend
known) results on the Boyd's indides, tensor, Sobolev embedding,
multiplicative and convolution operators in $ BSGLS $ spaces. }

The paper is organized as follows. In the next section we start
with an exemplary case just to give a feeling of what happens and
then present the main results on the so-called Boyd index of $
BSGLS \ - $ spaces, on tensors and multiplicative inequalities,
and the Sobolev embedding theorem and convolution inequalities,
each of these in a separate subsection. Further, the section
follows where we prove the statements, and in the last section we
discuss the sharpness of the obtained results.
\bigskip

\section{Results}

To get a flavor of the setting we work with, let, for instance, $
X =\mathbb R^n,$ $\sigma$ be a constant, and $\mu = \mu_{\sigma}$
be the measure on Borel subsets of $ X $ with density $
d\mu_{\sigma}/dx = |x|^{\sigma}.$ As usual, $x=
(x_1,x_2,\ldots,x_n)\in X$ so that $|x| = (x_1^2 + x_2^2 + \ldots
+ x_n^2)^{1/2}.$ let $L = L(z),$ $z \in (0,\infty) $ be a slowly
varying as $ z \to \infty $ positive continuous function, $ I(A) =
I(A,x) = 1$ for $x \in A$ and $I(A,x) = 0 $ otherwise. Let

\begin{eqnarray*}
f_L(x) &=& f(x; L,a,\alpha) = I(|x| > 1) \ |x|^{-1/a} \
[\log(|x|)]^{\alpha} \ L(\log |x|),\\
g_L(x) &=& g(x; L,b,\beta) = I(|x| < 1) \ |x|^{-1/b} \
[|\log(|x|)|]^{\beta} \ L( | \log |x| |),\\
A = A(a,n,\sigma) &=& a(n+\sigma) \ge 1, \ B = B(b,n,\sigma) = b(n
+ \sigma) \in (A, \infty),\\
\gamma &=& \alpha + 1/A, \ \delta = \beta + 1/B, \ p \in (A,B),
\end{eqnarray*}          and

\begin{eqnarray*}
\psi_L(p) = \psi_L(p; A,B; \gamma,\delta) \stackrel{def}{=}
(p-A)^{- \gamma} \ (B-p)^{-\delta} \ \max(L(A/(p - A)), L(B/(B -
p))).                                           \end{eqnarray*}
The function $ h_L(x) = h_L(x; a,b; \alpha, \beta) = f_L(x) +
g_L(x) $ belongs to the space $ G(\psi_L): $

\begin{eqnarray*}
h_L(\cdot) \in G(A,B; \gamma,\delta) \stackrel{def}{=}
G(\psi_L(\cdot)),                          \end{eqnarray*} and
this inclusion is exact in the sense that for $ p \in (A,B)$ there
holds  $|h_L|_p \asymp \psi_L(p),$ where here and in what follows
for $ p \in (A,B)$ the relation $g(p) \asymp h(p) $ denotes

\begin{eqnarray*}
0<\inf_{p\in (A,B)} f(p)/g(p) \le\sup_{p
\in(A,B)}f(p)/g(p)<\infty.                       \end{eqnarray*}
At the endpoints we need more in the case when $f(p) \to \infty.$
This may occur when either $p\to A+$ or $p\to B-$ or in both
cases. In detail, this means that in the case when $ \psi(A+0) =
\infty$ while $\psi(B-0) < \infty$ there holds

\begin{eqnarray*}
\lim_{p \to A + 0} \psi(p)/\nu(p) = 1;
\end{eqnarray*}
in the case when $ \psi(B-0) = \infty$ while $\psi(A+0) < \infty$
there holds

\begin{eqnarray*}\lim_{p \to B - 0} \psi(p)/\nu(p) = 1;
\end{eqnarray*}
and when in both cases  $ \psi(A+0) = \psi(B-0) = \infty$ there
holds

\begin{eqnarray*}
\lim_{p \to A + 0} \psi(p)/\nu(p) = \lim_{p \to B-0}
\psi(p)/\nu(p)= 1.
\end{eqnarray*}

  Denoting now $ \omega(n) = \pi^{n/2}/\Gamma(n/2 + 1)$ and
$\Omega(n) = n \omega(n) = 2 \pi^{n/2} / \Gamma(n/2), $ we let

\begin{eqnarray*}
R = R(\sigma,n) =  [(\sigma + n)/\Omega(n)]^{1/(\sigma + n)}, \
\sigma + n > 0,                                  \end{eqnarray*}
and let $ h = h(|x|) $ be a non-negative measurable function
vanishing for $ |x| \ge R(\sigma,n).$ For $u \ge e^2$ let

\begin{eqnarray*}
\mu_{\sigma} \{x: h(|x|) > u) \} = \min(1, \exp \left(- W(\log u))
\right),                                           \end{eqnarray*}
where $ W = W(z) $ is a twice differentiable strictly
convex for $ z \in [2, \infty) $ and strictly increasing function.
Denoting by

\begin{eqnarray*}
W^*(p) = \sup_{z > 2}(pz - W(z))                \end{eqnarray*}
the Young - Fenchel transform of the function $ W(\cdot), $ we
define the function

\begin{eqnarray*}\psi(p) = \exp \left( W^*(p)/p \right).\end{eqnarray*}

 It follows from the theory of Orlicz's spaces ( \cite[p.p. 22 -
27]{Os13}) that if for $ p \in [a, \infty) $ we have $|h|_p \asymp
\psi(p),$ then $ h(\cdot) \in G(\psi; a, \infty) $ and $  G(\psi;
1, \infty)$ coincides with some exponential Orlicz space. \par

We will restrict ourselves to the case $ p \in ( a(n+\sigma), b(n
+ \sigma)), \ a(n + \sigma) \ge 1, \ b(n+\sigma) < \infty,$ with
$p \to a(n + \sigma) + 0.$ Denoting $ A = a(n + \sigma),$ let

\begin{eqnarray*}
f(x) = f_L(x) = I(|x| > 1) \ |x|^{-1/a} \ (\log |x|)^{\gamma}
 \ L( \log |x|). \end{eqnarray*}
Using then multidimensional polar coordinates and well-known
properties of slowly varying functions (\cite[Ch.1, Sect.1.4 -
1.5]{Se15}), we obtain

\begin{eqnarray*}
 \Omega(n) \ ||f||_p^p &=& \int_1^{\infty} r^{- p/a + n + \sigma - 1} \
 (\log r)^{\gamma \ p} \ L^p(\log r) \ dr\\
 &=&(p/a - n - \sigma)^{-\gamma p - 1} \ \int_0^{\infty} e^{-z} \
z^{\gamma p} \ L^p(a z /(p - A)) \, dz\\
& \sim& (p/a - n - \sigma)^{- \gamma p - 1} \ L^p(a/(p - A))
\int_0^{\infty} e^{-z} \ z^{\gamma p} \, dz\\
& =& (p/a - n - \sigma)^{ - \gamma p - 1} \ L^p(a/(p - A)) \
\Gamma(1 + \gamma \ p).                         \end{eqnarray*}
Thus, we have for $ p \in (A,B)$

\begin{eqnarray*}
|f_L|_p \asymp (p - A)^{ - \gamma - 1/A} \ L(a/(p - A)).
\end{eqnarray*}

\subsection{Indices.}

 In this section we give an expression for
the so-called Boyd's (and other) indices of $ G(\psi,a,b) $ spaces
in the case of $X = [0, \infty) $ with usual Lebesque measure.
These indices play very important role in the theory of interpolation of
operators, in Fourier Analysis on r.i. spaces, etc. (see, e.g.,
\cite [p.p. 22 - 31, 192 - 204]{BS1}).

We recall the definitions. Given the family of (linear) operators
$ \{\sigma_s \} $ acting from some r.i. space $ G $ to $ G $ by the
following definition:

\begin{eqnarray*}
\sigma_s f(x) = f(x/s), \ s > 0, \ ||\sigma_s|| = ||\sigma_s||_{G
\to G}.                                          \end{eqnarray*}

 We have for arbitrary r.i. space $ G $ on the set $ X $ by the classical
definition ([1], chapter 2, [11], chapter 2)

\begin{eqnarray*}
\gamma_1(G) = \lim\limits_{s \to 0+} \log ||\sigma_s||/\log s;  \
&&\gamma_2(G) = \lim\limits_{s \to \infty} \log ||\sigma_s||/\log s;\\
\gamma^{(2)}(G) = \limsup\limits_{s \to 0+} \phi(G;2s)/\phi(G;s);
\ &&\gamma^{(1)}(G) = \liminf\limits_{s \to
0+}\phi(G;2s)/\phi(G;s).\end{eqnarray*}

\begin{theorem}\label{new1}  {\it There holds}

\begin{eqnarray*}
\gamma_1(G(\psi; a,b))= 1/b, \ \gamma_2(G(\psi; a,b))= 1/a, \ \psi
\in \Psi;                                        \end{eqnarray*}
{\it and if the space } $ X $ {\it is arbitrary and the measure} $ \mu $
{\it is diffuse}

\begin{eqnarray*}
\gamma^{(1)}( G(\psi; a,b) )= \gamma^{(2)}(G(\psi; a,b))= 2^{1/b}.
\end{eqnarray*} \end{theorem}

In a more general case of $ X = R^n$ with $\mu = \mu_{\sigma}$ and
$\sigma \ge 0,$ we analogously have

\begin{eqnarray*}
\gamma_1(G(\psi; a,b))= (n + \sigma)/b, \ \gamma_2(G(\psi; a,b))=
(n + \sigma)/a, \ \psi \in \Psi.
\end{eqnarray*}

\subsection{Tensor and multiplicative inequalities}

Let $ (X, \Sigma_1, \mu) $ and $ (Y, \Sigma_2,
 \nu) $ be two measurable spaces with $ \sigma$-finite measures
$ \mu$ and $\nu$ respectively. Let $ f = f(x) \in
G_X(\psi_1; a_1,b_1)$ and $g = g(y) \in G_Y(\psi_2; a_2,b_2),$
where $x \in X,$ $y \in Y,$ $\psi_1, \psi_2 \in \Psi,$ and let
 $a =  \max(a_1, a_2) < \min(b_1, b_2) = b.$ We set $\psi(p)
\stackrel{def}{=} \psi_1(p) \ \psi_2(p)$ for $p \in (a,b).$

 Let us consider the so-called {\it tensor product} of
 $ f,g: \ z(x,y) \stackrel{def}{=} f(x) \ g(y). $
Since both functions $ f$ and $g $ are independent of the space $
(X \times Y, \ \Sigma_1 \times \Sigma_2, \ \zeta),$ with $\zeta =
\mu \times \nu, $ we have:

\begin{lemma}\label{lem1} The following tensor inequality holds
\begin{eqnarray*}
||z||_{G(\psi; a,b)} \le ||f||_{G(\psi_1; a_1, b_1)}
||g||_{G(\psi_2; a_2, b_2)}.     \end{eqnarray*} \end{lemma}

This inequality is sharp, for example,  when $ \psi(p) = |f|_p$
and $\psi_2(p) = |g|_p. $ \par
  We consider now the so-called {\it multiplicative} inequality.
 Let  $ f  \in G(\psi_1; a_1, b_1)$ and $g \in G(\psi_2; a_2,
b_2),$ with $a_1,a_2 \ge 1$ and $1/b_1 + 1/b_2
> 1. $ We denote

\begin{eqnarray*}
A_1 = \max(1, a_1 a_2/(a_1 + a_2)),\quad B_1 = b_1 b_2/(b_1 +
b_2), \end{eqnarray*}             and

\begin{eqnarray*} \psi_3(r) =
\inf \{\psi_1(pr) \ \psi_2(qr); \ p,q: \ p,q > 1, 1/p + 1/q = 1
\}, \end{eqnarray*} with $r  \in (A_1, B_1).$

\begin{theorem}\label{th1}
There holds

\begin{eqnarray}\label{in1}
||f \, g||_{G(\psi_3; A_1, B_1)} \le ||f||_{G(\psi_1; a_1, b_1)} \
||g||_{G(\psi_2; a_2, b_2)}.
\end{eqnarray}

\end{theorem}

We mention that the sharpness of (\ref{in1}) up to multiplicative
constant  can be seen from
letting $ f = f_L \in G(\psi_1)$ and $g = g_L \in G(\psi_2)$ with
 $ \psi_1 = \psi_L(p; a_1,b_1; \alpha_1,\beta_1)$ and $\psi_2 =
\psi_L(p, a_2, b_2; \alpha_2, \beta_2). $ Namely, in the considered
case

\begin{eqnarray}\label{in1}
||f_L \, g_L||_{G(\psi_3; A_1, B_1)} \ge C(L, a_1, a_2, b_1, b_2) \
||f_L||_{G(\psi_1; a_1, b_1)} \ ||g_L||_{G(\psi_2; a_2, b_2)}.
\end{eqnarray}

We note, in addition, that if $ f \in G(\psi; a,b) $ and for $
\gamma = const  \in [a,b]$ we have $\psi_{\gamma}(p) =
\psi^{\gamma}(\gamma \, p), $ then $ g(x)
 = f^{\gamma}(x) \in G(\psi_{\gamma}; a/\gamma, b/\gamma) $ and

\begin{eqnarray*}
|| f^{\gamma}||_{G(\psi_{\gamma})} = ||f||^{\gamma}_ {G(\psi)}.
\end{eqnarray*}

\subsection{Sobolev embedding and convolution operators}

Let $ X $ be a convex domain in $ R^n, n \ge 2,$ with smooth
boundary, $ B $ be a projection operator on the $ m$-dimensional
smooth convex sub-manifold $ Y $ of $ X, \ m \le n, $ endowed with
the corresponding surface measure, and let $\psi = \psi(p; a,b) ,
\  1 \le a < b < n.$ We denote $ A_2 = \max(1, am/(n - a)),$ $B_2
= bm/(n - b),$ $\nu(q) = q^{1 - 1/n}\psi(qn/(q + m)),$ with $q \in
(A_2, B_2),$ and $u = u(x) \in C_1^0(X), $ i.e. $ u(\cdot) $ is
continuous differentiable and $ \lim_{|x| \to \infty} u(x) = 0.$

\begin{theorem}\label{th2} The following Sobolev type inequality holds

\begin{eqnarray}\label{in2}
||Bu||_{G(\nu; A_2, B_2)} \le C(X,Y; \psi) \, || \, |{\mbox {\rm
grad}} \, u| \, ||_{G(\psi; a,b)}. \end{eqnarray}
\end{theorem}

This result is supplied with the following interesting assertion,
opposite to that for the classical $L_p$ spaces.

\begin{theorem}\label{cor} In Theorem \ref{th2} the corresponding
embedding Sobolev operator is not compact. \end{theorem}

We now consider the (bilinear) generalized convolution
operator of the form

\begin{eqnarray*}
v(x)= (f*g)(x) = \int_X g (xy^{-1}) f(y) \ \mu(dy),
\end{eqnarray*}
where $ X $ is an unimodular Lie's group, $ \mu $ is its Haar's
measure. The unimodularity means, in particular, that $ \mu $ is
bide-side invariant. For the commutative group $ X,$ with standard
notation $ y^{-1} = - y,$ $xy^{-1} = x - y, $  this definition
coincides with the classical definition of convolution.\par
 Let $ f \in G(\psi_1; a_1, b_1)$ and $g \in G(\psi_2;
a_2, b_2)$ provided $1/a_1 + 1/a_2 > 1$ and $1/b_1 + 1/b_2 > 1. $
We denote

\begin{eqnarray*}
A_3 = a_1 a_2/(a_1 + a_2 - a_1 a_2), \ B_3 = b_1 b_2 /(b_1 + b_2 -
b_1 b_2),                                       \end{eqnarray*}
and for the values $ r \in (A_3, B_3) $ we define

\begin{eqnarray*}
\tau(r) =\inf \{\psi_1(p) \ \psi_2(q); \ p,q: \ p,q > 1, 1/p + 1/q
= 1 + 1/r \}. \end{eqnarray*}

\begin{theorem}\label{th3} There holds
\begin{eqnarray}\label{in3}
||f * g||_{G(\tau; A_3, B_3)} \le ||f||_{G(\psi_1; a_1, b_1)} \
||g||_{G(\psi_2; a_2, b_2)}.
\end{eqnarray}  \end{theorem}

\bigskip

\section{Proofs}

{\bf Proof of Theorem \ref{new1}.}

The last assertion of the theorem follows from the  explicit
expression for the fundamental function  $ \phi(\delta, G(\psi;
a,b)),$ with $\delta \in (0,\infty):$

\begin{eqnarray*}
 \phi (\delta; G) =  \sup \{||I(A)||_G, \ A \in \Sigma, \ \mu(A)
\le \delta \};                                  \end{eqnarray*}
see \cite{OS14}. The first assertion follows immediately from the
identity

\begin{eqnarray}\label{id}
|| \sigma_s || = \max \left( s^{1/a}, s^{1/b} \right), \ s > 0,
\psi  \in \Psi,
\end{eqnarray}
where

\begin{eqnarray*}
|| \sigma_s || \stackrel{def}{=} || \sigma_s||_{G(\psi) \to
G(\psi) }.
\end{eqnarray*}
It remains to prove (\ref{id}). The upper bound is obtained as
follows. \par

 Let $ f: f \in G(\psi), \ f \ne 0. $ We have

\begin{eqnarray*}
| \sigma_s f|^p_p = \int_0^{\infty} |f(x/s)|^p \ dx = s \
\int_0^{\infty} |f(y)|^p \ dy.
\end{eqnarray*}
It follows from

\begin{eqnarray*}
|\sigma_s f|_p &=& s^{1/p} |f|_p \le \max \left(s^{1/a}, s^{1/b}
\right) \ |f|_p\\
&\le& \max \left(s^{1/a}, s^{1/b} \right) \ \psi(p) \
||f||_{G(\psi)};
\end{eqnarray*}
that

\begin{eqnarray*}
|| \sigma_s || \le \max \left(s^{1/a}, s^{1/b} \right).
\end{eqnarray*}

For the lower bound, let $ g(\cdot) $ be a representation of
$ \psi: \ |g|_p=\psi(p),$ $p\in(a,b);$ then $||g||_{G(\psi)}=1$ and

\begin{eqnarray*}
|| \sigma_s || \ge || \sigma_s \ g||_{ G(\psi)} &=& \sup_{ p \in
(a,b)} [ \ |\sigma_s \ g|_p /\psi(p) \ ]\\
=\sup_{p \in (a,b) } [ s^{1/p} |g|_p /\psi(p)] &=& \sup_{p \in
(a,b) } s^{1/p} = \max \left(s^{1/a}, s^{1/b} \right).
\end{eqnarray*}

The proof is complete.    \hfill $\Box$

The proofs of our theorems \ref{th1}, \ref{th2} and \ref{th3} go
along similar lines and are strongly based on definitions and
preliminary matter given above.

{\bf Proof of Theorem \ref{th1}.} Let $ f \in G(\psi_1; a_1,b_1) $
and $ g \in G(\psi_2; a_2, b_2). $ By definition of these spaces,
this means

\begin{eqnarray*}
|f|_p &\le& \psi_1(p) \ \cdot ||f||_{G(\psi_1)}, \ p \in (a_1, b_1);\\
|g|_q &\le& \psi_2(q) \ \cdot ||g||_{G(\psi_1)}, \ q \in (a_2,
b_2).
\end{eqnarray*}
 It follows from H\"older's  inequality that for $ r \in (A_1, A_2) $
and $ p,q > 1,$ $1/p + 1/q = 1,$

\begin{eqnarray*}
|f \ g|_r \le |f|_{pr} \ |g|_{qr} \le \psi_1(pr) \ \psi_2(qr) \,
||f||_{G(\psi_1)} \ ||g||_{G(\psi_2)}.        \end{eqnarray*}

Minimizing the right-side over $ p$ and $q$ provided $p,q > 1$ and
$1/p + 1/q = 1,$ we obtain the desired assertion. \hfill $\Box$

{\bf Proof of Theorem \ref{th2}.} Here we will use the known
Sobolev inequality in the $ L_p $ spaces (see, e.g., \cite[Part 2,
Ch.11, Sect.4]{KA} or for newer and more extended versions
\cite{Tal17}) and \cite{You6}) which can be rewritten as

\begin{eqnarray*}
| Bu|_q \le C_1(X,Y) \, q^{1 - 1/n} \, | \, |{\mbox{\rm grad}} \,
u| \, |_p, \quad p = qn/(q+m). \end{eqnarray*}
  Let $  |{\mbox{\rm grad}} \ u| \in G(\psi; a,b)$ with $ b < n $ (the case $ b = n $
can be treated analogously); then we have for the values $ q \in
(A_2, B_2)$

\begin{eqnarray*}
| Bu|_q &\le& C_2(X,Y) || \, |{\mbox{\rm grad}} \, u| \,
||_{G(\psi)} \, \ q^{1 - 1/n} \ \psi(qn/(q + m))\\
&\le& C_2(X,Y) \, || \, |{\mbox{\rm grad}} \, u| \, ||_{G(\psi)}
\, \ \nu(q);
\end{eqnarray*} and
\begin{eqnarray*}
|| Bu||_{G(\nu)} \le C_2(X,Y) \, || \, |{\mbox{\rm grad}} \, u| \,
||_{G(\psi)},
\end{eqnarray*} This completes the proof.   \hfill $\Box$

{\bf Proof of Theorem \ref{th3}.} Assume that $ f \in G(\psi_1;
a_1,b_1) $ and $ g \in G(\psi_2; a_2, b_2). $ It follows from the
definition of these spaces that for $ p \in (a_1, b_1)$ and $q \in
(a_2, b_2) $

\begin{eqnarray*}
|f|_p \le ||f||_{G(\psi_1)}\, \psi_1(p), \quad |g|_q \le
||g||_{G(\psi_2)} \, \psi_2(q). \end{eqnarray*}
 Using the classical Young inequality (\cite{BL2}), we obtain

\begin{eqnarray*}
|f*g|_r \le C(p,q) \ |f|_p \,  |g|_q, \quad C(p,q) \le 1, \ 1 +
1/r = 1/p + 1/q,                              \end{eqnarray*}
where for $ n =  \dim \ X,$ $s = p/(p-1),$ $t = q/(q-1),$ and $ z =
r/(r-1)$ the constant $C(p,q)$ is given in \cite{BL2} in the case
 $ X = R^n $ in explicit form

\begin{eqnarray*}
C(p,q) = \left[p^{1/p} s^{-1/s} q^{1/q} t^{-1/t} r^{1/r} z^{-1/z}
\right]^{n/2}.\end{eqnarray*}   This yields

\begin{eqnarray*}
|f*g|_r \le \psi_1(p) \, \psi_2(q) \,||f||_{G(\psi_1)}
||g||_{G(\psi_2)}. \end{eqnarray*}
 The proof can now be completed as above, by
 minimizing over $ p $ and $ q, $ where $ p,q > 1$ and $1/p + 1/q =
1 + 1/r. $ \hfill $\Box$

{\bf Proof of Theorem \ref{cor}.} We introduce the Sobolev-Grand
Lebesque spaces $ W_1(\psi)$ with the (finite) norm of a function
$ u = u(x),$ defined on  $ X,$

\begin{eqnarray*}
||u||_{W_1(\psi)} = || \, |{\mbox{\rm grad}}\, u | \, ||_{G(\psi)}
+ || u ||_{G(\psi)}.
\end{eqnarray*} By this, the classical Sobolev embedding operator
$ S: \ W_1(\psi) \to G(\nu),$ $Su = Bu, $ is not compact, unlike
in the case of classical $ L_p $ spaces.\par
 This fact follows from the assertion that for the considered
function $ u = u(|x|) $ the family of  "small" shifts $
\{T_{\varepsilon} \, u(|x|) = u(\varepsilon + |x|) \},$ with
$\varepsilon \in (0, \varepsilon_0),$ $\varepsilon_0> 0,$ has a
positive distance in both $G(\psi)$ and $W_1(\psi)$ spaces: there
exists $C>0$ such that

\begin{eqnarray*}
||T_{\epsilon} u - T_{\delta} u||_{G(\psi)} \ge C
\end{eqnarray*}        and

\begin{eqnarray*}
||T_{\epsilon} u - T_{\delta} u||_{W_1(\psi)} \ge C,
\end{eqnarray*}
where $\epsilon, \delta \in (0, \epsilon_0),$ $\epsilon \ne
\delta.$

In order to prove the first assertion, (the second may be proved
analogously), we introduce certain
subspaces of $ G(\psi; a,b) $ space. Let us denote

\begin{eqnarray*}
GA(\psi)=GA(\psi; a,b) = \{f: \lim_{\delta \to 0+} \sup_{A: \mu(A)
\le \delta } ||f \ I(A) ||_{G(\psi)} = 0  \};      \end{eqnarray*}
and $GB(\psi) = GB(\psi; a,b)$ as the set of all $f$ such that for
all $\varepsilon> 0$ there exist $B \in (0,\infty)$ and $A \in
\Sigma,$ with $\mu(A) \le B,$ and there exists $g: \ X \to
{\mathbb R}$ such that $g(x) = g(x)I(A)$ and  $\sup_x |g(x)| < B$
and $||f - g|| < \varepsilon.$ Let also

\begin{eqnarray*}
G^0(\psi; a,b) = G^0(\psi) = \{f: \ \lim_{\psi(p) \to \infty}
|f|_p/\psi(p) = 0 \}.
\end{eqnarray*}
 The spaces $ GA(\psi), GB(\psi), G^0(\psi) $ are closed subspaces of
$ G(\psi). $ We assume here that there is $ h: X \to R,$ $|h|_p
\asymp \psi(p),$ $p \in (a,b).$

 It follows from the theory of r.i. spaces (\cite[Ch.1, p.p.22 - 28]{BS1})
  that if $ \psi \in \Psi,$ then

\begin{eqnarray*}
GA(\psi) = GB(\psi) = G^0(\psi) \ne G(\psi). \end{eqnarray*}
 Without loss of generality we can assume $ X = (0, 2 \pi],$ $\sigma = 0,$ and
define also  $x \pm y = x \pm y (\mbox{\rm mod} \, 2 \pi )$ for $
x,y \in X. $

 We take a function $ u = u(x) \in G(\psi) \setminus G^o(\psi), \ x \in X. $
Then (see \cite[Ch.3, p.p.192 - 198]{BS1})

\begin{eqnarray*}
\inf_{\epsilon \ne \delta} ||T_{\varepsilon}u -
T_{\delta}u||_{G(\psi)} = \inf_{\varepsilon
\ne\delta}||T_{\varepsilon - \delta}u-u||_{G(\psi)}> 0.
\end{eqnarray*}
The proof is complete.   \hfill $\Box$

\bigskip

\section{Sharpness}

We will discuss either the sharpness or lack of that for obtained
results. Since the sharpness of Theorems \ref{new1}, \ref{th1} and
\ref{th2} is obtained readily, we did not postpone it and gave
immediately after the formulations.

Let us demonstrate (briefly) the sharpness of Theorem \ref{th3} by
 considering only the case $ n = 1$ and $\sigma = 0. $
 Let both $ \gamma_1, \gamma_2 \ge 0,$ and define

\begin{eqnarray*}
f(x) &=& I(0 < x < 1) \ x^{ - 1/b_1} \ |\log x|^{\gamma_1}, \ b_1,
b_2 > 1,\\ g(x) &=& I(0 < x < 1) \ x^{- 1/b_2} \ |\log
x|^{\gamma_2}, \ h(t) = (f*g)(t). \end{eqnarray*}
 It suffices to consider the case $ t \in (0, 1/4), $ since on
$ [1/4, 2] $ the function $ h = h(t) $ is bounded and for $ t \in
( - \infty,0) \cup (2, \infty)$ this function vanishes.

For $ t \to 0+$ we get

\begin{eqnarray*}
h(t) &=& \int_0^t x^{- 1/b_1} |\log x|^{\gamma_1} \ (t - x)^{-
1/b_2} \ |\log(t - x)|^{\gamma_2} \ dx\\ &=& t^{1 - 1/b_1 - 1/b_2}
\int_0^1 y^{- 1/b_1} \ (1 - y)^{-1/b_2} \, |\log t + \log
y|^{\gamma_1} \ |\log t + \log(1 - y) |^{\gamma_2} dy\\
& \sim & t^{1 - 1/b_1 - 1/b_2} \ |\log t|^{\gamma_1 + \gamma_2} \
\int_0^1 y^{ - 1/b_1} \ (1 - y)^{-1/b_2} \, dy\\
& =& B(1 - 1/b_1, 1 - 1/b_2) \ t^{1 - 1/b_1 - 1/b_2} \ |\log
t|^{\gamma_1 + \gamma_2}; \end{eqnarray*}
 here  $ B(\cdot, \cdot) $ stands for the beta-function.
Taking then the $|\cdot|_p$ norm, we obtain

\begin{eqnarray*}
|f*g|_p \sim C ( B_3 - p) ^{ - \gamma_1 - \gamma_2 - 1/b_1 - 1/b_2
+ 1}, \ p \in [1, B_3), \end{eqnarray*}

The situation is more complicated with (\ref{in2}) and
(\ref{in3}). Roughly speaking, we can present examples when the
inequalities are achieved but for the same examples the actual
bounds are better. In other words, the sharpness of these
inequalities is an open problem.

We first analyze (\ref{in2}). Let $ \psi(p)= (p - a)^{- \alpha} \
(b - p)^{ - \beta},$ $\sigma = 0,$ $1 \le a < b < n. $ Consider
the function $ u = u(|x|),$ with $X = R^n$ and $Y = R^m,$ for
which $|{\mbox{\rm grad}} \ u| \in G(a,b; \alpha,\beta), $ i.e.,
for $p \to a+0$ and $p \to b - 0$ we have $| \ |{\mbox{\rm grad}}
\ u| \ |_p \sim (p - a)^{ - \alpha} \ (b - p)^{-\beta},$ then it
follows from the inequality (\ref{in2}) that

\begin{eqnarray*}
|Bu|_p \le C (p - A_2)^{- \alpha} \ (B_2 - p)^{- \beta}, \ p \in
(A_2, B_2).                                      \end{eqnarray*}
However, in the considered case for the same range of $p$

\begin{eqnarray*}
|Bu|_p \sim C (p - A_2)^{ - \alpha + 1/n} \, (B_2 - p)^{ - \beta +
1/n}.\end{eqnarray*}

In the same way, analogous examples may be constructed in the
cases when either $ b_1 = \infty $ or $ b_2 = \infty. $ \par

 Therefore, the {\it bounds} $ A_2$ and $B_2 $ are, in general,
exact, but between the {\it exponents} obtained  $ - \alpha$ and
$- \beta $ on one side and $ - \alpha + 1/n$ and $- \beta + 1/n $
in the example on the other side there is the $ 1/n$ "gap".\par

 But since the dimension $ n $ may be sufficiently great, we can
conclude that  in general case the assertion (3) of theorem 4 is
exact.\par

We give a similar example for (\ref{in3}). Let $ X = R,$ $\sigma =
0,$ $\gamma_1, \gamma_2 \ge 0,$ $1 \le a_1, a_2 < b_1, b_2 <
\infty,$ and $1/b_1 + 1/b_1> 1. $ Considering the functions

\begin{eqnarray*}
f(x) &=& I(0 < |x| < 1) \ |x|^{-1/b_1} \ |\log |x| \ |^{\gamma_1},\\
g(x) &=& I(0 < |x| < 1) \ |x|^{-1/b_2} \ |\log |x| \ |^{\gamma_2},
\end{eqnarray*}
and $h = f*g,$ we then have

\begin{eqnarray*}
f \in G(1,b_1; 0, \gamma_1 + 1/b_1), \ g \in G(1,b_2; 0, \gamma_2
+ 1/b_2).                                        \end{eqnarray*}
It follows from (\ref{in3}) that for $ p \in [1, B_3) $

\begin{eqnarray*}
|h|_p \le C (B_3 - p) ^{ - \gamma_1 - \gamma_2 - 1/b_1 - 1/b_2}.
\end{eqnarray*}
In fact

\begin{eqnarray*}
|h|_p \sim C (B_3 - p) ^{- \gamma_1 - \gamma_2 - 1/b_1 - 1/b_2 +
1}.                                         \end{eqnarray*}
Analogously, if $ h = f*g, $ where

\begin{eqnarray*}
f(x) &=& I( |x| > 1) \ |x|^{ - 1/a_1} \ | \log |x| \ |^{\gamma_1},\\
g(x) &=& I( |x| > 1) \ |x|^{ - 1/a_2} \ | \log |x| \ |^{\gamma_2},
\end{eqnarray*}
then
\begin{eqnarray*}
f \in G(a_1, b_1; \gamma_1 + 1/a_1, 0), \ g \in G(a_2, b_2;
\gamma_2 + 1/b_2, 0).                              \end{eqnarray*}
It follows from (\ref{in3}) that

\begin{eqnarray*}
|h|_p \le C (p - A_3) ^{ - \gamma_1 - \gamma_2 - 1/a_1 - 1/a_2},
\end{eqnarray*}
but in fact

\begin{eqnarray*}\ |h|_p \sim C (p - A_3) ^{ - \gamma_1 - \gamma_2 - 1/a_1 -
1/a_2 + 1}.\end{eqnarray*}

 Therefore, like above the {\it bounds} $ A_3$ and $B_3 $ are in this case
exact, but between the {\it exponents} there is a 1 "gap".\par
 Note that in the considered case the "gap" does not tends to zero as
 $ n \to \infty, $ opposite to the previous (convolution) case.\par

 Finding sharp estimates in these two cases is an interesting open
problem.


\begin{thebibliography}{99}

\bibitem{BS1}
C. Bennet and R. Sharpley, {\it Interpolation of operators.}
Orlando, Academic Press Inc., 1988.

\bibitem{BL2}
H.J. Brascamp and E.H. Lieb, {\it Best constants in Young's
inequality, its converse and its generalization to more than three
functions.} Journ. Funct. Anal., {\bf 20}(1976), 151--173.

\bibitem{CM3}
M. Carro and J. Martin, {\it Extrapolation theory for the real
interpolation method.} Collect. Math. {\bf 33}(2002), 163--186.

\bibitem{Fio4}
A. Fiorenza. {\it Duality and reflexivity in grand Lebesgue
spaces.} Collect. Math. {\bf 51}(2000), 131--148.

\bibitem{FK5}
A. Fiorenza and G.E. Karadzhov, {\it Grand and small Lebesgue
spaces and their analogs.} Consiglio Nationale Delle Ricerche,
Instituto per le Applicazioni del Calcoto Mauro Picine", Sezione
di Napoli, Rapporto tecnico 272/03(2005).

\bibitem{IS7}
T. Iwaniec and C. Sbordone, {\it On the integrability of the
Jacobian under minimal hypotheses.} Arch. Rat.Mech. Anal., {\bf
119}(1992), 129--143.

\bibitem{IKO8}
T. Iwaniec, P. Koskela and J. Onninen, {\it Mapping of finite
distortion: Monotonicity and Continuity.} Invent. Math. {\bf
144}(2001), 507--531.

\bibitem{JM9}
B. Jawerth and M. Milman, {\it Extrapolation theory with
applications.} Mem. Amer. Math. Soc. {\bf 440}(1991).

\bibitem{KA}
L.V. Kantorovich and G.P. Akilov, {\it Functional analysis. Second
edition}. Pergamon Press, Oxford-Elmsford, N.Y. 1982.

\bibitem{KM10}
G.E. Karadzhov and M. Milman, {\it Extrapolation theory: new
Results and applications.} J. Approx. Theory, {\bf 113}(2005),
38--99.

\bibitem{KO11}
Yu.V. Kozatchenko and E.I. Ostrovsky, {\it Banach spaces of random
variables of subgaussian type.} Theory Probab. Math. Stat., Kiev,
1985,  42--56.

\bibitem{KPS12}
S.G. Krein, Yu. Petunin and E.M. Semenov, {\it Interpolation of
Linear operators.} New York, AMS, 1982.

\bibitem{LedTal18}
M. Ledoux and M. Talagrand  (1991) {\it Probability in
 Banach Spaces.} Springer, Berlin, 1991.

\bibitem{Os13}
E.I. Ostrovsky, {\it Exponential Estimations for Random Fields.}
Moscow - Obninsk, OINPE, 1999 (Russian).

\bibitem{OS14}
E. Ostrovsky and L.Sirota, {\it Some new rearrangement invariant
spaces: theory and applications.} Electronoc publications:
arXiv:math.FA/0605732 v1, 29 May 2006.

\bibitem{You6}
Y.J. Park, {\it Logarithmic Sobolev trace inequality.} Proc. Am.
Math. Soc. {\bf 132}(2004), 2075--2083.

\bibitem{Se15}
E. Seneta, {\it Regularly Varying Functions}. Springer Verlag,
Berlin - Heidelberg - New York, 1976.


\bibitem{Tal17}
G. Talenti, {\it Inequalities in Rearrangement Invariant
Functional Spaces.} Nonlinear analysis, Function Spaces and
Applications. Prometheus, Prague, {\bf 5}(1995), 177-230.

%\bibitem{You6}
%Young Ja Park, {\it Sobolev Trace Inequalities.} Electronic
%Publications, arXiv: math.CA/0107065 v1 9 Jul. 2001.




\end{thebibliography}
\end{document}